\documentclass[11pt,twoside]{amsart}
\usepackage[T1]{fontenc}
\usepackage[utf8]{inputenc}
\usepackage[a4paper]{geometry}
\usepackage{amsmath,amssymb,amsfonts,amsthm,mathtools,mathrsfs}
\usepackage{graphicx}
\usepackage{xcolor}
\usepackage{cite}
\usepackage{tikz-cd}
\usepackage[bb=stix]{mathalpha}
\usepackage{fancyhdr}
\usepackage{hyperref}
\hypersetup{colorlinks=true,linkcolor=blue,citecolor=blue,urlcolor=blue,pdftitle={On Functoriality of Baum--Bott Residues},pdfauthor={Mauricio Correa and Tatsuo Suwa},pdfsubject={Baum--Bott residues and functoriality},pdfkeywords={Baum--Bott residues, holomorphic foliations, functoriality, localized characteristic classes, rationality}}

\newcommand{\Q}{{\mathbb{Q}}}

\newcommand{\C}{{\mathbb{C}}}
\newcommand{\Z}{{\mathbb{Z}}}

\def\C{\mathbb C}

\renewcommand{\dim}{\mathrm{dim}}

\newcommand{\F}{\mathscr F}

\newcommand{\G}{\mathscr G}

\def\Z{\mathbb Z}

\newcommand{\Sing}{{\rm Sing}}

\newcommand{\Res}{{\rm Res}}

\newcommand{\cod}{{\rm cod}}

\newcommand{\R}{{\mathbb{R}}}

\newcommand{\mfd}{{manifold}}
\newcommand{\ra}{{\rightarrow}}
\newcommand{\hra}{{\hookrightarrow}}
\newcommand*{\Cdot}{\raisebox{-0.4ex}{\scalebox{1.8}{$\hspace{.3mm}\cdot\hspace{.3mm}$}}}
\newcommand*{\mCdot}{\raisebox{-0.4ex}{\scalebox{1.8}{$\hspace{.3mm}\cdot\hspace{-.2mm}$}}}

\newtheorem{lema}{Lemma}[section]
\newtheorem{cor}[lema]{Corollary}
\newtheorem{teo}[lema]{Theorem}

\theoremstyle{definition}
\newtheorem{remark}[lema]{Remark}

\newtheorem{defi}[lema]{Definition}
\newtheorem{exe}[lema]{Example}

\begin{document}

\pagestyle{fancy}

\title{On functoriality of Baum--Bott residues}

\author{Maur\'icio Corr\^ea}
\address{Maur\'icio Corr\^ea  \\
Universit\`a degli Studi di Bari,
Via E. Orabona 4, I-70125, Bari, Italy
}
\email[M. Corr\^ea]{mauricio.barros@uniba.it,mauriciomatufmg@gmail.com }

\author{ Tatsuo Suwa}
\address{ Tatsuo Suwa\\ Department of Mathematics  \\Hokkaido University \\
Sapporo 060-0810, Japan}
\email[Tatsuo Suwa]
{tsuwa@sci.hokudai.ac.jp }



\begin{abstract}
We prove a functoriality formula for Baum--Bott residues under holomorphic pull-back. Let $f:V\to X$ be a holomorphic map generically transverse to a singular holomorphic foliation $\F$ on $X$. Under a compatibility condition on the tangent and normal sheaves, the localized Baum--Bott class of $f^*\F$ is the pull-back of that of $\F$, and the corresponding residue is obtained by localized intersection.  We also discuss the universal relative class of Baum and Bott and its relation with the rationality conjecture. If the characteristic polynomial has rational coefficients, rational residues remain rational under the pull-backs considered here. For proper surjective generically finite maps of equal dimension we obtain, in addition, a push-forward formula; hence rationality is equivalent before and after pull-back, in particular for bimeromorphic modifications satisfying the hypotheses of the theorem.
\end{abstract}

\maketitle

\hyphenation{sin-gu-lari-ties}
\hyphenation{re-si-dues}

\tableofcontents


\section{Introduction}\label{secintro}

Baum and Bott \cite{BB} defined residues associated with singular holomorphic foliations. Their explicit determination is already delicate in higher dimension. For a one-dimensional foliation, the residue at an isolated singularity is expressed by a Grothendieck residue. For a compact connected component of the singular set of the expected dimension, Baum and Bott obtained the following transverse description.

Let $X$ be a complex manifold of dimension $n$ and $\F$ a singular holomorphic foliation of dimension $k$ on $X$. Let $S$ be a connected component of the singular set $\Sing(\F)$. For $p\in S$, choose holomorphic vector fields $v_1,\dots,v_s$ defined on an open neighbourhood $U_p$ of $p$ in $X$ so that, for every $x\in U_p$, the germs at $x$ of $v_1,\dots,v_s$ are in $T\F_x$ and span $T\F_x$ as an $\mathcal O_x$-module, where $T\F$ denotes the tangent sheaf of $\F$. Let $\mathbf T_p(\F)$ be the linear subspace of $T_pX$ spanned by $v_1(p),\dots,v_s(p)$ and set
\[
S^{(i)}=\{\,p\in S\mid \dim\mathbf T_p(\F)\leq k-i\,\},\qquad i=1,\dots,k.
\]
Then there is a filtration
\[
S=S^{(1)}\supset S^{(2)}\supset\cdots\supset S^{(k)}.
\]

Assume now that $S$ has the ``expected dimension'', i.e., $\dim S=k-1$, and let $S_1,\dots,S_s$ be its irreducible components of dimension $k-1$. Baum and Bott proved in \cite[Theorem~3]{BB}, under the assumptions
\[
\dim S=k-1\qquad\text{and}\qquad \dim S^{(2)}<k-1,
\]
that, for a symmetric homogeneous polynomial $\phi$ of degree $n-k+1$, the residue $\Res_{\phi}(\F,S)$ is a linear combination of the classes $[S_1],\dots,[S_s]$, with coefficients given by the corresponding transverse residues at nonsingular points of the $S_i$ away from $S^{(2)}$.

We recall the corresponding transverse formula. Let $S'$ be an irreducible component of $S$ of dimension $k-1$, let $p\in S'\setminus S^{(2)}$ be nonsingular, and consider an inclusion $i:\mathbb D^{n-k+1}\to U_p$ of a sufficiently small disc of dimension $n-k+1$, centred at $0\in\C^{n-k+1}$ and transverse to $S'$ at $i(0)=p$. Let $i^*\F$ denote the foliation on $\mathbb D^{n-k+1}$ induced from $\F$ by $i$. It is one-dimensional and has $\{0\}$ as an isolated singularity. Thus $\Res_{\phi}(i^*\F,0)\in\C$ is represented by a Grothendieck residue at $0$. The above result of Baum and Bott gives
\begin{equation}\label{BBmax}
\Res_{\phi}(i^*\F,0)=\left(\mathbb D^{n-k+1}\Cdot\Res_{\phi}(\F,S')\right)_0,
\end{equation}
where the right-hand side is the intersection product of $\mathbb D^{n-k+1}$ with $\Res_{\phi}(\F,S')$, localized at $0$; see Section~\ref{sec:Baum-Bott} below.

In \cite{CorreaLourenco2019}, the first named author and Louren\c co showed that the residue for singularities of the expected dimension can be computed in the same manner without the above generic condition, but still under a compactness assumption. In \cite{Vishik}, formula \eqref{BBmax} was obtained under the assumption that the tangent sheaf of the foliation is locally free. These results extend the transverse computation in different directions.

Formula \eqref{BBmax} is a transverse pull-back formula. Baum and Bott observed that their local residue classes are ``functorial in an appropriate sense'' and related this functoriality to the homotopy and homology of the classifying spaces $B\Gamma_q^{\C}$ of complex foliations \cite[p.~288]{BB}. At the end of the proof of the transverse formula they noted that the argument ``really verifies a very special case of the functoriality'' of the residue \cite[p.~328]{BB}.

Section~12 of \cite{BB} gives the corresponding homotopy-theoretic formulation. Put $q=n-k$. There is a classifying space $B\Gamma_q^{\C}$ and a natural map
\begin{equation}\label{eq:nu-classifying}
\nu:B\Gamma_q^{\C}\longrightarrow BGL(q),
\end{equation}
which assigns to a regular complex foliation its normal bundle. Composing with $BGL(q)\to BGL$, one obtains a map to the stable classifying space. Since only its homotopy type is relevant, Baum and Bott regard this map as an inclusion, or equivalently work with its mapping cylinder and the pair $(BGL,B\Gamma_q^{\C})$.

If $\phi$ is a symmetric homogeneous polynomial of degree $d>q$, the Bott vanishing theorem says that the corresponding universal characteristic class restricts trivially to $B\Gamma_q^{\C}$. Baum and Bott interpret the construction of their $Z$-sequences as assigning to $\phi$ a definite relative lifting
\begin{equation}\label{eq:universal-class}
\widetilde\phi\in H^{2d}(BGL,B\Gamma_q^{\C};\C).
\end{equation}
This relative class is the one used in the localization of $\phi$ at the singular set.

The same relative classes occur in the rigidity and rationality questions. In the critical degree $d=q+1$, the transverse formulas and the examples of Section~11 of \cite{BB} may retain continuous transverse parameters. For $q+1<d<n$, Baum and Bott state that no general explicit formula is known, prove their rigidity theorem, and formulate the rationality conjecture: if $\phi$ has rational coefficients, the corresponding residue should be rational. Examples~3--5 of \cite{BB} show the contrast: critical-degree residues may vary with the defining parameters, while higher-degree residues are rigid and, in the examples that can be computed explicitly, rational when $\phi$ is rational. The functoriality theorem below also applies in this higher-degree range.

The rationality conjecture has also been studied by localization methods. Suwa \cite{Suwa1984} proved rationality for a class of residues of foliations of complete-intersection type in terms of Chern classes and local Chern classes. Sert\"oz \cite{Sertoz1989} used the Nash modification and the Grassmann graph construction to reduce the problem, under suitable hypotheses, from coherent sheaves to vector bundles. Brasselet and Suwa \cite{BrasseletSuwa2000} compared Baum--Bott and Nash residues and obtained further cases of the conjecture. These results concern Nash-type modifications, whereas the present paper concerns holomorphic pull-back.

At the universal level one should distinguish the existence of a rational lifting from the particular lifting selected by the Baum--Bott construction. If $\phi$ has rational coefficients, Bott vanishing holds with rational coefficients and the long exact sequence of the pair gives a relative lifting
\[
\widetilde\phi_{\Q}\in H^{2d}(BGL,B\Gamma_q^{\C};\Q).
\]
This does not settle the rationality conjecture, since the Baum--Bott construction singles out the definite complex lifting in \eqref{eq:universal-class}. The difference between this distinguished lifting and the complexification of a rational lifting lies in the image of
\[
H^{2d-1}(B\Gamma_q^{\C};\C)\longrightarrow H^{2d}(BGL,B\Gamma_q^{\C};\C).
\]
The rationality problem therefore concerns the rational structure of the distinguished localized class, not only the existence of a rational relative lifting. Its image in relative cohomology with coefficients in $\C/\Q$ would have to vanish in order to obtain a universal form of the rationality statement. Without a realization or detection result for geometric classifying maps, this universal condition is stronger than the Baum--Bott conjecture itself.

Let $\F$ be singular on $X$ and let $S$ be a connected component of $\Sing(\F)$. After restricting to a neighbourhood $U$ of $S$ containing no other singular points, the regular foliation on $U\setminus S$ gives a classifying map to $B\Gamma_q^{\C}$, while a real analytic locally free resolution of the normal sheaf gives the corresponding stable classifying map to $BGL$. On the regular part, the exact sequence of vector bundles supplies the homotopy between these two maps. Thus one obtains, up to homotopy, a map of pairs
\[
c_{\F,S}:(U,U\setminus S)\longrightarrow(BGL,B\Gamma_q^{\C}),
\]
and the Baum--Bott construction gives
\begin{equation}\label{eq:universal-localization}
\phi_S(N\F,\F)=c_{\F,S}^{*}(\widetilde\phi).
\end{equation}
The Baum--Bott residue is obtained from this localized class by Alexander duality. Under this duality, formula~\eqref{BBmax} is the transverse case of the functoriality considered by Baum and Bott.

Let $V$ be a complex manifold and let $f:V\to X$ be a holomorphic map generically transverse to $\F$, and write $\G=f^*\F$. On the regular locus, the normal bundle and its basic connection behave naturally under transverse pull-back. At the singular set, however, a locally free resolution of $N\F$ need not remain exact after pull-back, and $N\G$ need not coincide with $f^*N\F$. We give conditions under which both difficulties are avoided and the localized class pulls back.

If $S$ is a connected component of $\Sing(\F)$ and $Z=f^{-1}(S)$, we give conditions under which
\begin{equation}\label{eq:intro-localized}
\phi_Z(N\G,\G)=f^*\phi_S(N\F,\F)
\qquad\text{in }H^{2d}(V,V\setminus Z;\C).
\end{equation}
By Alexander duality, this yields the corresponding functoriality formula
\begin{equation}\label{eq:intro-functoriality}
\Res_{\phi}(\G,Z)=\left(V\mCdot_f\Res_{\phi}(\F,S)\right)_Z.
\end{equation}
The intersection product on the right is the localized intersection product with respect to $f$ recalled in Section~\ref{sec:Baum-Bott}. No compactness assumption is involved in this local statement.

The issue is the behaviour of the normal sheaf and of its resolutions under pull-back. The following compatibility condition ensures that the localized class on $V$ is represented by the pull-back of the data on $X$.

Our main result is the following; see Sections~\ref{sec.fol} and~\ref{sec:Baum-Bott} for the notation and terminology.

\begin{teo}\label{Main-Theorem}
Let $X$ be a complex manifold of dimension $n$ and $\F$ a singular holomorphic foliation of dimension $k$ on $X$. Let $V$ be a complex manifold of dimension $n'\geq n-k+1$ and let $f:V\to X$ be a holomorphic map generically transverse to $\F$. Let $S$ be a connected component of $\Sing(\F)$, set $Z=f^{-1}(S)$, and denote by $\G=f^*\F$ the induced foliation. Assume that
\begin{enumerate}
\item[(i)] there is an open set $U\supset S$, with $(U\setminus S)\cap\Sing(\F)=\emptyset$, such that $T\F$ admits a resolution by real analytic vector bundles on $U$;
\item[(ii)] there is an open set $W\supset Z$ such that $(W\setminus Z)\cap\Sing(\G)=\emptyset$;
\item[(iii)] on a neighbourhood of $Z$, the canonical sequence on the regular locus extends to an exact sequence
\begin{equation}\label{eq:normal-compatibility-intro}
0\longrightarrow f^*T\F\longrightarrow f^*TX\longrightarrow N\G\longrightarrow0.
\end{equation}
\end{enumerate}
Then, for every symmetric homogeneous polynomial $\phi$ of degree $d\geq n-k+1$, the residues $\Res_{\phi}(\F,S)$ and $\Res_{\phi}(\G,Z)$ are defined and
\begin{equation}\label{formula}
\Res_{\phi}(\G,Z)=\left(V\mCdot_f\Res_{\phi}(\F,S)\right)_Z
\qquad\text{in }H^{\rm BM}_{2n'-2d}(Z;\C).
\end{equation}
\end{teo}

In the universal description above, the content of the theorem is that the classifying data associated with $\G$ can be chosen compatibly with the pull-back of those associated with $\F$. Indeed, condition \rm{(iii)}, together with the Tor vanishing proved in Section~\ref{sec:Main-Theorem}, shows that the real analytic resolution defining the localized class of $\F$ pulls back to a genuine resolution of $N\G$. On the regular complement, the canonical isomorphism $N\G\simeq f^*N\F$ identifies the classifying map of the induced foliation with the pull-back of that of $\F$. Thus, on sufficiently small neighbourhoods of $S$ and $Z$, the corresponding maps of pairs are compatible with $f$, and \eqref{eq:universal-localization} gives
\[
c_{\G,Z}^{*}(\widetilde\phi)=f^*c_{\F,S}^{*}(\widetilde\phi),
\]
which is \eqref{eq:intro-localized}. Thus, under the hypotheses of Theorem~\ref{Main-Theorem}, the universal relative Baum--Bott class is compatible with pull-back. Under Alexander duality this gives the commutative diagram
\[
\begin{tikzcd}
H^{2d}(U,U\setminus S;\C)\arrow[r,"A_{U,S}"]\arrow[d,"f^*"'] & H^{\rm BM}_{2n-2d}(S;\C)\arrow[d,"(V\mCdot_f\ )_Z"]\\
H^{2d}(f^{-1}(U),f^{-1}(U)\setminus Z;\C)\arrow[r,"A_{f^{-1}(U),Z}"'] & H^{\rm BM}_{2n'-2d}(Z;\C),
\end{tikzcd}
\]
where the localized intersection product is taken with respect to the restriction of $f$ to $f^{-1}(U)$. By excision this is the same localized class as in the statement of the theorem. Thus formula \eqref{formula} is the Borel--Moore homological form of the naturality of the universal localized class.

We record the corresponding consequence for the rationality problem.

\begin{cor}\label{cor:rationality-pullback}
Under the hypotheses of Theorem~\ref{Main-Theorem}, suppose that $\phi$ has rational coefficients. If
\[
\Res_{\phi}(\F,S)\in\operatorname{im}\bigl(H^{\rm BM}_{2n-2d}(S;\Q)\longrightarrow H^{\rm BM}_{2n-2d}(S;\C)\bigr),
\]
then
\[
\Res_{\phi}(\G,Z)\in\operatorname{im}\bigl(H^{\rm BM}_{2n'-2d}(Z;\Q)\longrightarrow H^{\rm BM}_{2n'-2d}(Z;\C)\bigr).
\]
In particular, in the range $q+1<d<n$, the rationality property conjectured by Baum and Bott is preserved by the pull-backs covered by Theorem~\ref{Main-Theorem}.
\end{cor}

\begin{proof}
The localized intersection operation in \eqref{formula} is the Borel--Moore homological counterpart, under Alexander duality, of the pull-back in relative cohomology. Both pull-back and Alexander duality are defined with rational coefficients. Hence the localized intersection operation sends rational classes to rational classes, and the assertion follows from \eqref{formula}.
\end{proof}

Equivalently, let
\[
\rho_{\phi}(\F,S)
=
[\Res_{\phi}(\F,S)]
\in
\frac{H^{\rm BM}_{2n-2d}(S;\C)}{\operatorname{im}H^{\rm BM}_{2n-2d}(S;\Q)}
\]
be the rationality defect of the residue, and define $\rho_{\phi}(\G,Z)$ similarly. Theorem~\ref{Main-Theorem} gives
\begin{equation}\label{eq:rationality-defect}
\rho_{\phi}(\G,Z)
=
\left(V\mCdot_f\rho_{\phi}(\F,S)\right)_Z.
\end{equation}
Equation~\eqref{eq:rationality-defect} shows that the rationality defect is functorial. This does not prove the Baum--Bott rationality conjecture; it shows only that rationality is preserved by the pull-backs covered by the theorem. For proper surjective generically finite maps of equal dimension, Corollary~\ref{cor:proper-degree} gives the converse implication, so the original residue is rational if and only if the pulled-back residue is rational. For lower-dimensional pull-backs, a family detecting the quotient containing $\rho_{\phi}(\F,S)$ would give a reduction of the rationality question. Example~\ref{ex:finite-etale-rationality} gives an explicit finite \'etale example in the higher-degree range.

\begin{remark}\label{rem:intro-main}
{\bf 1.} Condition \rm{(i)} is local near $S$. It is satisfied, for example, if $T\F$ is locally free near $S$; if $S$ is contained in a relatively compact open set on which a real analytic locally free resolution can be chosen, cf.~\cite[Section~6]{BB}; and for germs, after shrinking a representative and tensoring a holomorphic free resolution with the sheaf of real analytic functions.
\smallskip

\noindent
{\bf 2.} Condition \rm{(ii)} only requires that $Z=f^{-1}(S)$ be isolated, in a neighbourhood of $Z$, from the other components of $\Sing(\G)$.
\smallskip

\noindent
{\bf 3.} Suppose in addition that $\operatorname{codim}_V Z\geq2$. Then condition \rm{(iii)} is satisfied if $f^*T\F$ and $N\G$ are reflexive on a neighbourhood of $Z$; see Lemma~\ref{lem:normal-sequence}. In particular, it is enough to assume that $f^*T\F$ is reflexive and that
$
\mathcal Ext^1(\Omega_V^1/N\G^\vee,\mathcal O_V)=0,
$
since the latter condition implies that $N\G$ is reflexive by Lemma~\ref{lem:normal-reflexive}. The codimension assumption is automatic when $Z\subset\Sing(\G)$. Hence the Ext condition is sufficient for the normal compatibility required in \rm{(iii)}.
\smallskip

\noindent
{\bf 4.} None of $X$, $V$ or $S$ is required to be compact. The construction is localized: the characteristic class lies in relative cohomology $H^{2d}(U,U\setminus S;\C)$ and Alexander duality identifies it with a class in the Borel--Moore homology of $S$. When $S$ is compact, the residue agrees with the usual Baum--Bott residue. Compactness is needed in the original global residue formula in order to recover the global characteristic class from the sum of the local residues, but not for the local class itself.
\end{remark}

If $V$ is a complex submanifold of $X$ and $f=i:V\hookrightarrow X$ is the inclusion, formula \eqref{formula} becomes
\begin{equation}\label{eq:intro-inclusion}
\Res_{\phi}(i^*\F,V\cap S)=\left(V\Cdot\Res_{\phi}(\F,S)\right)_{V\cap S}.
\end{equation}
If $\operatorname{codim}_V(V\cap S)\geq2$, condition \rm{(iii)} follows whenever $i^*T\F$ and $N(i^*\F)$ are reflexive. A sufficient condition is that $T\F$ be locally free near $S$ and that
$
\mathcal Ext^1(\Omega_V^1/N(i^*\F)^\vee,\mathcal O_V)=0.
$
Indeed, Remark~\ref{rem:intro-main}.3 then applies. The same reflexivity criterion holds for arbitrary holomorphic maps under the corresponding codimension assumption.

Theorem~\ref{Main-Theorem} contains the transverse formula of Baum and Bott as a particular case. Indeed, under the assumptions $\dim S=k-1$ and $\dim S^{(2)}<k-1$, their local structure theorem shows that, at a point $p\in S_{\rm reg}\setminus S^{(2)}$, the foliation singularity is locally the pull-back, by a submersion, of an isolated zero of a holomorphic vector field. For a sufficiently small normal disc $i:\mathbb D^{n-k+1}\hookrightarrow X$, the normal sheaf of the induced one-dimensional foliation is identified with the restriction of $N\F$, and condition \rm{(iii)} is satisfied. Formula \eqref{formula} therefore reduces to \eqref{BBmax}. In contrast, Theorem~\ref{Main-Theorem} does not require $S$ to have the expected dimension, does not involve the condition $\dim S^{(2)}<k-1$, is not restricted to the critical degree $n-k+1$, and applies to holomorphic maps rather than only to transverse normal discs.

The result of \cite{CorreaLourenco2019} removes the genericity condition in the expected-dimensional computation under global hypotheses. Theorem~\ref{Main-Theorem}, instead, concerns functoriality under holomorphic maps and does not require compactness. When condition \rm{(iii)} is satisfied it applies, in particular, to the corresponding transverse sections. The locally free case is also related to the functorial formula considered by Vishik \cite{Vishik}.

The proof uses the real analytic resolutions of Baum and Bott. Condition \rm{(iii)} implies the vanishing of the first Tor sheaf associated with the pull-back of $N\F$. After factoring $f$ through its graph, Lichtenbaum's rigidity theorem over regular local rings gives the vanishing of all higher Tor sheaves. Hence the resolution defining the localized class of $\F$ pulls back to a resolution of $N\G$. Equality \eqref{eq:intro-localized} follows from the naturality of Chern--Weil forms, and \eqref{formula} from Alexander duality.

The paper is organised as follows. In Section~\ref{sec.fol} we recall singular holomorphic foliations, transverse pull-backs and the normal sheaf. In Section~\ref{sec:Baum-Bott} we recall the localized intersection product and Baum--Bott residues for possibly non-compact components. Theorem~\ref{Main-Theorem} is proved in Section~\ref{sec:Main-Theorem}, and Section~\ref{sec:consequences} records some immediate consequences.

\subsection*{Acknowledgments}
MC is grateful to Hokkaido University for its hospitality. He is partially supported by the Universit\`a degli Studi di Bari and by PRIN 2022MWPMAB, ``Interactions between Geometric Structures and Function Theories'', and is a member of INdAM--GNSAGA. He was partially supported by CNPq grants 202374/2018-1 and 400821/2016-8, and by FAPEMIG grants APQ-02674-21, APQ-00798-18 and APQ-00056-20. MC would like to thank Javier Gargiulo Acea, Olivier Thom, Alan Muniz, Jos\'e Seade, Pablo Perrella and Sebasti\'an Velazquez for useful discussions and corrections. TS is partially supported by JSPS grants 20K03572 and 23K25772..

\section{Singular holomorphic foliations}\label{sec.fol}
In this section, $X$ denotes a complex manifold of dimension $n$ and $\mathcal O_X$ the sheaf of holomorphic functions on $X$. We do not distinguish holomorphic vector bundles and locally free $\mathcal O_X$-modules.

A singular holomorphic \textit{distribution} $\F$ of dimension $k$ on $X$ is given by an exact sequence
\begin{equation}\label{eq:foliation}
\begin{tikzcd}
0\arrow[r] & T\F\arrow[r,"\iota"] & TX\arrow[r,"\pi"] & N\F\arrow[r] & 0,
\end{tikzcd}
\end{equation}
where $N\F$ is torsion-free. If $T\F$ is involutive, $[T\F,T\F]\subset T\F$, we say that $\F$ is a singular holomorphic \textit{foliation}. The sheaves $T\F$ and $N\F$ are the tangent and normal sheaves of $\F$, respectively. The codimension of $\F$ is $q=n-k$.

The singular set is
\[
\Sing(\F):=\Sing(N\F)=\bigcup_{i=1}^{n-1}\operatorname{Supp}\mathcal Ext^i(N\F,\mathcal O_X).
\]
Since $N\F$ is torsion-free, $\Sing(\F)$ has codimension at least two, and $T\F$ is reflexive; cf.~\cite{Kobayashi}. 
The conormal sheaf is $N\F^\vee$. Taking the maximal exterior power of the dual of $\pi$ gives a twisted $q$-form
\[
\omega_\F\in H^0(X,\Omega_X^q\otimes\det(N\F))
\]
whose zero set is $\Sing(\F)$. On $X\setminus\Sing(\F)$ it is locally decomposable and defines the regular foliation in the usual way.

\begin{defi}\label{def:generically-transverse}
Let $V$ be a complex manifold of dimension greater than $q$. A holomorphic map $f:V\to X$ is \textit{generically transverse} to $\F$ if $f^*\omega_\F$ is not identically zero and
\[
df_y(T_yV)+T\F_{f(y)}=T_{f(y)}X
\]
whenever $f^*\omega_\F(y)\ne0$.
\end{defi}

In this case $f^*\omega_\F$ defines a singular holomorphic foliation $\G=f^*\F$ on $V$, of the same codimension $q$, with
\[
T\G=\left(\{v\in TV\mid\imath_v(f^*\omega_\F)=0\}\right)^{\vee\vee}.
\]
Let $V^0\subset V$ be an open set on which $\G$ is regular, $\F$ is regular along $f(V^0)$, and $f$ is transverse to $\F$. The morphism
\[
TV^0\longrightarrow f^*N\F|_{V^0},\qquad v\longmapsto f^*\pi(df(v)),
\]
is surjective and has kernel $T\G|_{V^0}$. Hence there is a canonical isomorphism
\begin{equation}\label{eq:normal-regular}
\chi:N\G|_{V^0}\overset{\sim}{\longrightarrow}f^*N\F|_{V^0}.
\end{equation}

Recall that the Bott partial connection of a regular foliation is natural under transverse pull-back. Thus, if $\nabla$ is a basic connection for $N\F$ on an open set on which $\F$ is regular, the connection on $N\G$ corresponding to $f^*\nabla$ by \eqref{eq:normal-regular} is basic. This may be checked in local Frobenius coordinates: if the leaves of $\F$ are the fibres of a submersion $h$, then the leaves of $\G$ are the fibres of $h\circ f$, and the two Bott partial connections correspond under $N\G\simeq f^*N\F$.

\begin{lema}\label{lem:normal-reflexive}
Let $\G$ be a singular holomorphic foliation on a complex manifold $V$. If
$
\mathcal Ext^1(\Omega_V^1/N\G^\vee,\mathcal O_V)=0,
$
then $N\G$ is reflexive.
\end{lema}

\begin{proof}
Dualizing the exact sequence
\[
0\longrightarrow N\G^\vee\longrightarrow\Omega_V^1\longrightarrow\Omega_V^1/N\G^\vee\longrightarrow0
\]
gives
\[
0\longrightarrow(\Omega_V^1/N\G^\vee)^\vee\longrightarrow TV\longrightarrow N\G^{\vee\vee}\longrightarrow0.
\]
Since $T\G$ is reflexive and agrees with $(\Omega_V^1/N\G^\vee)^\vee$ outside $\Sing(\G)$, whose codimension is at least two, the two sheaves are isomorphic. Comparing with $0\to T\G\to TV\to N\G\to0$ gives $N\G\simeq N\G^{\vee\vee}$.
\end{proof}

\section{Baum--Bott residues for non-compact components}\label{sec:Baum-Bott}

The Baum--Bott residue was originally defined for a compact connected component of the singular
set of a foliation and the residue formula, which says that the sum of the residues gives the
global characteristic class of the normal sheaf of the foliation, is stated for a foliation on a compact manifold, see \cite[Theorem 2]{BB}. In this section, we remove these compactness assumptions. For details of the material
in this section, we refer to  \cite{Suwa1998,Suwa2024} and the references therein. 

Let $X$ be a $C^\infty$ manifold of real dimension $m$. Recall that, for any closed sub\mfd\ $V$ of $X$, there exists a $C^\infty$ triangulation of $X$ compatible with $V$. 

Let $K_0$ be a $C^\infty$ triangulation of $X$, let $K$ be its barycentric subdivision, and let $K'$ be the barycentric subdivision of $K$. Thus $K'$ is the second barycentric subdivision of $K_0$. The use of the second subdivision ensures that the star, relative to $K'$, of a $K_0$-subcomplex $L$ is homotopy equivalent to its polyhedron $|L|$.
For a $p$-simplex $\mathbb{s}$ of $K$, denote by $\mathbb{s}^*$ the $(m-p)$-cell dual to $\mathbb{s}$, namely the union of the $(m-p)$-simplices of $K'$ meeting $\mathbb{s}$ at its barycentre $b_{\mathbb{s}}$. Let $K^*$ denote the resulting dual cell decomposition of $X$. We choose orientations for the simplices and dual cells. Singular homology and cohomology may then be represented by the corresponding simplicial and cellular complexes.

\subsection{Dualities} 
We denote by $H^p_{K^*}(X;\mathbb Z)$ the $p$-th cohomology of the cochain complex of finite $K^*$-cochains. It is canonically isomorphic to the singular cohomology $H^p(X;\mathbb Z)$. For homology we use locally finite $K$-chains; their homology is canonically the Borel--Moore homology $H_p^{\rm BM}(X;\mathbb Z)$.

Assume that $X$ is oriented. Poincar\'e duality gives the isomorphism (cf.~\cite{Br1,Suwa2024})
\[
\begin{tikzcd}
 P_X:H^p(X; \mathbb{Z})\simeq  H^p_{K^*}(X;\mathbb{Z}) \arrow[r,"\sim"] & H_{m-p}^{\rm BM}(X; \mathbb{Z}).
\end{tikzcd}
\]
At the chain level it sends a \(p\)-cochain $u$ of \(K^*\) to
\[
\sum_{  \mathbb{s}} \langle   \mathbb{s}^*, u\rangle \mathbb{s}, 
\]
where $\mathbb{s}$ runs through the $(m-p)$-simplices of $X$. Equivalently, it is the left cap product with the fundamental class of $X$:
\[
P_X([u]) = [X]  \frown [u].
\]
 

Let $S$ be a closed subset of $X$ and suppose that $K_0$ is compatible with $S$, so that $S$ is a $K_0$-subcomplex. 
Recall that the star \( \mathbb{S}_{K'}(S) \) of \( S \) in \( K' \) is the union of simplices of \( K' \) intersecting with \( S \), i.e., the union of cells of \( K^* \) intersecting with \( S \). Let
$
\mathbb{O}_{K'}(S) = \mathbb{S}_{K'}(S) \setminus \partial \mathbb{S}_{K'}(S)
$
denote the open star. 
Then there is a proper deformation retraction
$
\mathbb{S}_{K'}(S) \to S
$
and a deformation retraction
$
\mathbb{O}_{K'}(S) \to S.
$
Hence there are natural isomorphisms
\[
H^p_{K^*}(X, X \setminus \mathbb{O}_{K'}(S);\Z)\simeq H^p(X, X \setminus \mathbb{O}_{K'}(S);\Z) \simeq H^p(X, X \setminus S;\Z).
\]
In this situation, Alexander duality gives the isomorphism
\[
\begin{tikzcd}
A_{X, S}:H^p(X, X \setminus S;\Z)\simeq H^p_{K^*}(X, X \setminus \mathbb{O}_{K'}(S);\Z)  \arrow[r,"\sim"] & 
H_{m-p}^{\rm BM}(S;\Z).
\end{tikzcd}
\]
At the chain level it sends a \(p\)-cochain $u$ of \(K^*\) vanishing on $X\setminus\mathbb{O}_{K'}(S)$ to
\[
\sum_{  \mathbb{s}} \langle   \mathbb{s}^*, u\rangle \mathbb{s}, 
\]
where $\mathbb{s}$ runs through the $(m-p)$-simplices in $S$.



The following diagram is commutative; coefficients are omitted, and may be taken in $\Z$ or $\C$:
\[
\begin{tikzcd}
 H^p(X, X \setminus S) \arrow[d, "A_{X,S}"']  \arrow[r, "j^*"]  &  
 H^p(X)  \arrow[d, "P_X"] \\
[.2cm]  H_{m-p}^{\rm BM}(S) \arrow[r, "i_*"] &  
  H_{m-p}^{\rm BM}(X)
\end{tikzcd}
\]
where \(j:(X, \emptyset) \to (X, X \setminus S) \) and \(i:S \to X \)
denote the inclusions.
For more details, we refer the reader to \cite[Section 4.2]{Suwa2024}.

\subsection{\v{C}ech-de Rham theory} 
Let $X$ be a $C^\infty$ manifold of dimension $m$, and let $A^p(X)$ denote the complex vector space of complex-valued $C^\infty$ $p$-forms on $X$.
Consider the de Rham complex of \( X \), denoted by \( (A^\bullet(X), d) \), and the associated de Rham cohomology \( H_d^p(X) \). 

Let \( \mathcal{U} = \{ U_\alpha \}_{\alpha \in I} \) be an open covering of \( X \).
The \textit{\v{C}ech-de Rham cohomology} \( H^p_D(\mathcal{U}) \) is the \v{C}ech hypercohomology of the 
de~Rham complex on \( \mathcal{U} \).
For a covering $\mathcal U=\{U_0,U_1\}$ by two open sets, set $U_{01}=U_0\cap U_1$.
We set
\[
A^p(\mathcal{U}) = A^p(U_0) \oplus A^p(U_1) \oplus A^{p-1}(U_{01}).
\]
An element $\xi\in A^p(\mathcal U)$ is a triple $\xi=(\xi_0,\xi_1,\xi_{01})$, where $\xi_i$ is a $p$-form on $U_i$, $i=0,1$, and $\xi_{01}$ is a $(p-1)$-form on $U_{01}$. The differential $D:A^p(\mathcal U)\to A^{p+1}(\mathcal U)$ is
\[
D\xi = (d\xi_0, d\xi_1, \xi_1 - \xi_0 - d\xi_{01}).
\]

The morphism $A^p(X)\to A^p(\mathcal{U})$ given by $\omega\mapsto (\omega|_{U_0},\omega|_{U_1},0)$
induces an isomorphism (\cite[Theorem 7.1]{Suwa2024})
\[
\begin{tikzcd}
H_d^p(X) \arrow[r,"\sim"] & H_D^p(\mathcal{U}).
\end{tikzcd}
\]

\subsection{Relative \v{C}ech-de Rham cohomology}
Let \( X \) be a $C^\infty$ manifold of dimension \( m \) and \( S \) a closed set in \( X \). Let \( U_0 = X \setminus S \) and \( U_1 \) an open neighbourhood of \( S \). We then consider the covering \( \mathcal{U} = \{ U_0, U_1 \} \) of \( X \).  

Define
\[
\begin{aligned}
A^p(\mathcal U,U_0)
&:=\{\,\xi=(\xi_0,\xi_1,\xi_{01})\in A^p(\mathcal U)\mid \xi_0=0\,\}\\
&=A^p(U_1)\oplus A^{p-1}(U_{01}).
\end{aligned}
\]
Then $(A^\bullet(\mathcal U,U_0),D)$ is a subcomplex of $(A^\bullet(\mathcal U),D)$. The differential is given by
\[
D(\xi_1, \xi_{01}) = (d\xi_1, \xi_1 - d\xi_{01}).
\]
\begin{defi}
The \( p \)-th relative \v{C}ech-de Rham cohomology \( H_D^p(\mathcal{U}, U_0) \) of the pair \( (\mathcal{U}, U_0) \) is the \( p \)-th cohomology of \( (A^\bullet(\mathcal{U}, U_0), D) \).
\end{defi}

By \cite[Corollary 7.8]{Suwa2024} we have:

\begin{teo} 
If $S$ is a subcomplex of a triangulation of $X$, there is a canonical isomorphism
\[
\begin{tikzcd}
H_D^p(\mathcal{U}, U_0) \arrow[r,"\sim"] & H^p(X, X \setminus S;\C).
\end{tikzcd}
\]
\end{teo}

In terms of the dual cell decomposition, this isomorphism is described as follows.
Assume that $S$ is a subcomplex of a triangulation $K_0$ of $X$, and let $K$, $K'$ and $K^*$ be as above. We also take a ``honeycomb system $\{R_0,R_1\}$ adapted to  
$\mathcal{U}$, $K'$ and $S$". In our situation, 
we may set \( R_1 = \mathbb{S}_{K''}(S) \), the star of \( S \) in the barycentric subdivision \( K'' \) of \( K' \), and \( R_0 = X \setminus \text{Int}\,R_1 \). Set \( R_{01} = \partial R_0 \), which is equal to \( -\partial R_1 \), the boundary of \( R_1 \) with opposite orientation. 
Note that, for each cell \( \mathbb{s}^* \) of \( K^* \), the intersections \(\mathbb{s}^*\cap R_1\) and \(\mathbb{s}^*\cap R_{01}\) are $K''$-chains.

With these choices, the isomorphism is induced by the map sending a relative \v{C}ech-de Rham $p$-cochain $\xi=(\xi_1,\xi_{01})$ to the $p$-cochain of $K^*$ given by
\[
\mathbb{s}^*\mapsto 
 \int_{\mathbb{s}^* \cap R_1} \xi_1 + \int_{\mathbb{s}^* \cap R_{01}} \xi_{01}.
\]

If $X$ is oriented, Alexander duality with $\C$-coefficients takes the following form.

\begin{teo}\label{thAlexC}  
If $S$ is a subcomplex of a triangulation $K_0$ of $X$, there is a canonical isomorphism
\[
\begin{tikzcd}
A_{X,S}: H^p(X, X \setminus S;\mathbb{C})\simeq 
H_D^p(\mathcal{U}, U_0) \arrow[r,"\sim"] & H_{m-p}^{BM}( S;\mathbb{C}), 
\end{tikzcd}
\]
which assigns, to a relative class \( [\xi] = [(\xi_1, \xi_{01})] \),  the class of the cycle
\[
\sum_{\mathbb{s} \subset S} \left( \int_{\mathbb{s}^* \cap R_1} \xi_1 + \int_{\mathbb{s}^* \cap R_{01}} \xi_{01} \right) \mathbb{s},
\]
where \( \mathbb{s} \) runs through all the oriented \( (m-p) \)-simplices of \( K \) in \( S \).
\end{teo}

\subsection{Localized intersection product with a map} 
Let $X$ be an oriented 
$C^\infty$ manifold of dimension $m$. 
We take a triangulation $K_0$ of $X$ and let $K$, $K'$ and $K^*$ be as before. 
For a $p$-simplex ${\mathbb{t}}$ of $K'$ we denote by $\vartheta({\mathbb{t}})$ the $p$-cochain  
dual to ${\mathbb{t}}$, i.e., for every 
$p$-simplex ${\mathbb{t}}'$ of $K'$,
\[
\langle{\mathbb{t}'},\vartheta({\mathbb{t}})\rangle=\begin{cases}1&\quad\text{if}\ \ {\mathbb{t}'}={\mathbb{t}},\\
0&\quad\text{otherwise}.
\end{cases}
\]


Let ${\mathbb{s}}_{1}$ be an oriented $(m-r)$-simplex and  ${\mathbb{s}}_{2}$ an oriented $s$-simplex  of $K$.
We denote by $\frown_{r}$ and $\frown_{l}$ the right and left cap products of chains and cochains of 
$K'$, see \cite[Remark B.11]{Suwa2024}. Let $X_{K'}$ also denote the fundamental cycle of
$X$ in $K'$, i.e., the sum of the simplices of $K'$ in $X$.

The \textit{intersection product} ${\mathbb{s}}_{1}^{*}\Cdot {\mathbb{s}}_{2}$ is an 
$(r+s-m)$-chain of $K'$ defined by (cf.~\cite[Definition 4.3]{Suwa2024})
\[
{\mathbb{s}}_{1}^{*}\Cdot {\mathbb{s}}_{2}=(\vartheta({\mathbb{t}}_{1})\frown_{r}X_{K'})\frown_{l}\vartheta({\mathbb{t}}_{2}),
\]
where ${\mathbb{t}}_{1}$ and ${\mathbb{t}}_{2}$ are  an $(m-r)$-simplex and   an $(m-s)$-simplex  of $K'$, respectively, such that ${\mathbb{t}}_{1}\subset {\mathbb{s}}_{1}$ and ${\mathbb{t}}_{2}\subset {\mathbb{s}}_{2}^{*}$.

The following properties hold:
\begin{enumerate}
\item the definition
does not depend on the 
choice of ${\mathbb{t}}_{1}$ or ${\mathbb{t}}_{2}$,
\item the support of ${\mathbb{s}}_{1}^{*}\Cdot {\mathbb{s}}_{2}$ is the set-theoretic
intersection of  ${\mathbb{s}}_{1}^{*}$ and ${\mathbb{s}}_{2}$,
\item for any simplex ${\mathbb{s}}$ of $K$, we have ${\mathbb{s}}^*\Cdot {\mathbb{s}}=b_{\mathbb{s}}$, the barycenter of ${\mathbb{s}}$ (with multiplicity $1$).
\end{enumerate}

In general, we have the boundary formula
 \[
 \partial({\mathbb{s}}_{1}^{*}\Cdot {\mathbb{s}}_{2})= (-1)^{m-s}(\partial{\mathbb{s}}_{1}^{*})\Cdot {\mathbb{s}}_{2}+ {\mathbb{s}}_{1}^{*}\Cdot \partial{\mathbb{s}}_{2},
 \]
and hence induces the intersection product in homology
\[
\begin{tikzcd}
H_{r}^{\rm BM}(X;\Z)\times H_{s}^{\rm BM}(X;\Z)\arrow[r,"(\ \Cdot\ )"] & H^{\rm BM}_{r+s-m}(X;\Z).
\end{tikzcd}
\]
Under Poincar\'e duality it corresponds to the cup product in cohomology. 
More generally, let $S_1$ and $S_2$ be subcomplexes of $X$ and set $S=S_1\cap S_2$. There is a localized intersection product
\[
\begin{tikzcd}
H_{r}^{\rm BM}(S_{1};\Z)\times  H_{s}^{\rm BM}(S_{2};\Z)\arrow[r,"(\ \Cdot\ )_S"]& H_{r+s-m}^{\rm BM}(S;\Z).
\end{tikzcd}
\]

We also have the cup product
\[
\begin{tikzcd}
H^{p}(X,X\setminus S_1;\Z)\times H^{q}(X,X\setminus S_2;\Z)\arrow[r,"\smile"]& H^{p+q}(X,X\setminus S;\Z).
\end{tikzcd}
\]

The following diagram is commutative:

\[
\begin{tikzcd}
H^{p}(X,X\setminus S_1)\times H^{q}(X,X\setminus S_2) \arrow[d,"A_1\times A_2"']  \arrow[r,"\smile"]  &  H^{p+q}(X, X \setminus S) \arrow[d, "A" ]     \\
[.3cm] H_{r}^{\rm BM}(S_1)\times H_{s}^{\rm BM}(S_2) \arrow[r, "(\ \Cdot\ )_S"] &   H_{r+s-m}^{\rm BM}(S),
\end{tikzcd}
\]
where $p+r=m$, $q+s=m$ and $A$, $A_1$ and $A_2$ denote the Alexander isomorphisms for $S$, $S_1$ and $S_2$. The coefficients $\Z$ are omitted; the same diagram holds with $\C$-coefficients.

Let $V\subset X$ be a closed sub\mfd\ of dimension $m'$, and write $i:V\hra X$ for the inclusion.
By the tubular neighbourhood theorem, we may identify a neighbourhood of $V$ in $X$ with a neighbourhood of the zero section in the normal bundle of $V$. We assume that $V$ is oriented 
  and that the normal bundle of $V$ is orientable and oriented so that the orientation of the fibres, followed by the orientation of $V$, gives the orientation of $X$.
We may assume that $V$ is a $K_{0}$-subcomplex of $X$.

Intersection with $V$ is described by the following commutative diagrams (cf.~\cite[Proposition 4.13]{Suwa2024}):
\[
\begin{tikzcd}
H^{p}(X) \arrow[d,"i^*"']  \arrow[r,"P_X"]  &  H_{m-p}^{\rm BM}(X) \arrow[d, "(V\Cdot\ )_V" ]     \\
[.2cm] H^{p}(V) \arrow[r, "P_V"] &   H_{m'-p}^{\rm BM}(V),
\end{tikzcd}
\qquad
\begin{tikzcd}
H^{p}(X,X\setminus S) \arrow[d,"i^*"']  \arrow[r,"A_{X,S}"]  &  H_{m-p}^{\rm BM}(S) \arrow[d, "(V\Cdot\ )_Z" ]     \\
[.2cm] H^{p}(V,V\setminus Z) \arrow[r, "A_{V,Z}"] &   H_{m'-p}^{\rm BM}(Z),
\end{tikzcd}
\]
where $Z=V\cap S$.

\begin{figure}[htbp]
    \centering
    \includegraphics[width=0.62\textwidth]{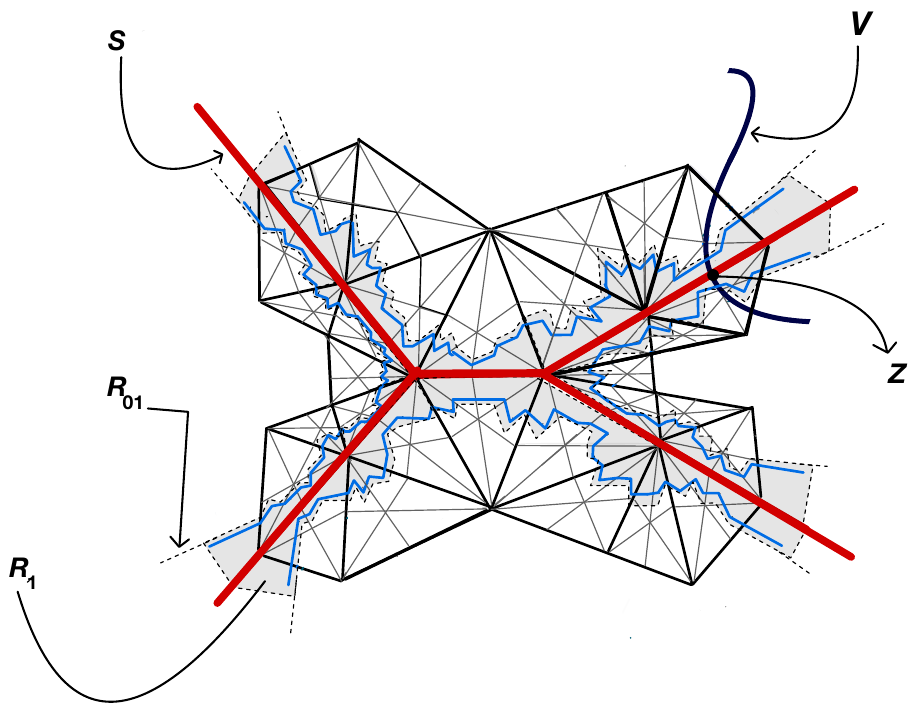}
    \caption{The honeycomb neighbourhood of $S$ and its localized intersection with $V$. The intersection $Z=V\cap S$ is shown together with $R_1$ and its boundary $R_{01}$.}
    \label{fig:honeycomb-intersection}
\end{figure}

The local configuration is illustrated in Figure~\ref{fig:honeycomb-intersection}. The region $R_1$ is a honeycomb neighbourhood of $S$, with $R_{01}=\partial R_0=-\partial R_1$, and the localized intersection is supported on $Z=V\cap S$.

Let $V$ be an oriented $C^\infty$ manifold of dimension $m'$ and let $f:V\to X$ be a $C^\infty$ map.
For a subcomplex $S$ of $X$, set $Z=f^{-1}(S)$ and assume that $Z$ is a subcomplex of a triangulation of $V$.

We define the \textit{intersection product}  \( (V \mCdot_f\ \ )_V \) and the \textit{localized intersection product}  \( (V \mCdot_f \ \ )_Z \) so that the following diagram is commutative (omitting the coefficient $\Z$, which could be $\C$):

\begin{equation}\label{lipm}
\begin{tikzcd}[row sep=1.3em, column sep=1.5em]
& H^{p}(X) \arrow[rr, "P_X"] \arrow[dd, "f^*"']  & & H_{m-p}^{\rm BM}(X) \arrow[dd, "(V \mCdot_f\ )_V"]  \\
H^{p}(X,X\setminus S) \arrow[ur, "j^*"]\arrow[rr, crossing over, "A_{X,S}" pos=0.7] \arrow[dd, "f^*"'] & & H_{m-p}^{\rm BM}(S) \arrow[ur, "i_*"] \arrow[dd, crossing over, "(V \mCdot_f\ )_Z" pos=0.3]   \\
& H^{p}(V) \arrow[rr, "P_V" pos=0.3]    & & H_{m'-p}^{\rm BM}(V) \\
H^{p}(V,V \setminus Z) \arrow[ur, "j^*"] \arrow[rr,  "A_{V,Z}"']   & & H_{m'-p}^{\rm BM}(Z) \arrow[ur, "i_*"]
\end{tikzcd}
\end{equation}

If \( V \) is a submanifold of \( X \) and \( f = i: V \hookrightarrow X \) is the inclusion map, then \( (V \mCdot_i\ )_V =(V \Cdot\ )_V\) and \( (V \mCdot_i\ )_Z =(V \Cdot\ )_Z\).

\begin{remark} In Fulton \cite{Fu}, similar notions of intersection products are defined, in the algebraic category, using normal cones.
\end{remark}

Let \( X \) and \( V \) be \( C^\infty \) manifolds, and let \( S \) and \( Z \) be subcomplexes of \( X \) and \( V \), respectively, with respect to triangulations of \( X \) and \( V \). 
Consider a \( C^\infty \) map \( f: (V, V \setminus Z)  \to (X, X \setminus S)\). We set \( U_0 = X \setminus S \) and \( W_0 = V \setminus Z \). Let \( U_1 \) and \( W_1 \) be neighbourhoods of \( S \) and \( Z \), respectively, such that \( f(W_1) \subset U_1 \). We then have the coverings \( \mathcal{U} = \{U_0, U_1\} \) and \( \mathcal{W} = \{W_0, W_1\} \) of \( X \) and \( V \), respectively. The map $f$ induces pull-backs
$ f^* : H_{D}^p(\mathcal{U}, U_0) \to H_{D}^p(\mathcal{W}, W_0) $
and
$f^* : H^p(X, X \setminus S; \mathbb{C}) \to H^p(V, V \setminus Z; \mathbb{C})$.
They fit into the commutative diagram
\[
\begin{tikzcd}
 H_{D}^p(\mathcal{U}, U_0) \arrow[d,"f^*"']  \arrow[r,"\sim"]  &  H^p(X, X \setminus S; \mathbb{C}) \arrow[d, "f^*" ]     \\
[.2cm]  H_{D}^p(\mathcal{W}, W_0) \arrow[r, "\sim"] &    H^p(V, V \setminus Z; \mathbb{C})
\end{tikzcd}
\]

\subsection{Baum--Bott residues}\label{ssBB}
Let $X$ be a complex manifold of dimension $n$. Note that, by the Cauchy--Riemann equations, the real tangent bundle $T_{\R}X$ of $X$ is canonically isomorphic to the holomorphic tangent bundle $TX$ and that the complexification $T_{\R}^cX=\C\otimes T_{\R}X$ is canonically isomorphic to the direct sum $TX\oplus\overline{T}X$.
For an open set $U$ in $X$, we denote by \( A^0(U) \) the \( \C \)-algebra of \( C^\infty \) functions on \( U \). For a complex vector bundle \( E \) on \( X \), we set 
$A^p(U;E)=C^{\infty}(U;\bigwedge^p(T_{\R}^cX)^\vee\otimes E)$, the $A^0(U)$-module of $C^{\infty}$ $p$-forms with coefficients in $E$. 
 
Recall that a connection for \( E \) is a \( \C \)-linear map
\[
\begin{tikzcd}
\nabla:A^0(X; E) \arrow[r] & A^1(X; E)
\end{tikzcd}
\]
satisfying the Leibniz rule:
\[
\nabla(f s) = df \otimes s + f \nabla(s) \qquad \text{for}\ \ f \in A^0(X),\ s \in A^0(X; E).
\]

Let \( \F \) be a singular holomorphic foliation of dimension \( k \) on \( X \). Note that, if $(N\F)_x$ is $\mathcal{O}_{X,x}$-free, so is $(T\F)_x$. Thus, setting 
$X_0=X\setminus \Sing(\F)$,
 there is an involutive subbundle $TF$ of $TX_0$  such that $T\F|_{X_0}=\mathcal{O}_{X_0}(TF)$. If we set $NF=T{X_0}/TF$, then we have $N\F|_{X_0}=\mathcal{O}_{X_0}(NF)$. There is an action of $TF$ on $NF$:
\[
\begin{tikzcd}
 TF\times NF\arrow[r] & NF\qquad\text{given by}\ \  (u,\pi(v))\mapsto \pi([u,v]),
 \end{tikzcd}
\]
where $\pi:TX_0\ra NF$ denotes the canonical surjection.

A connection for $NF$ on $X_0$ is said to be \textit{basic} if it is compatible with the action, i.e.,
\begin{equation}\label{basic}
\nabla(\pi(v))(u)=\pi([u,v])\qquad\text{for}\ \ v\in A^0(X_0;TX_0),\ u\in A^0(X_0;TF),
\end{equation}
and is of type
$(1,0)$. The Bott vanishing theorem says that, if $\nabla$ is a basic connection for $NF$, then
\[
\phi(\nabla)=0
\]
for every symmetric homogeneous polynomial \( \phi \) of degree \( d \ge n - k+1 \).

Denote by \( \mathcal{A}_X \) the sheaf of  real analytic functions on \( X \). 
We do not distinguish real analytic vector bundles and locally free $\mathcal{A}_X$-modules.

Let $S$ be a connected component of $\Sing(\F)$ and 
suppose that there is a neighbourhood $U$ of $S$, disjoint from the other components, such that   \( N\F|_U \) admits a resolution by real analytic vector bundles:
$$
\begin{tikzcd}
0 \arrow[r] & E_r \arrow[r] & E_{r-1} \arrow[r] & \cdots \arrow[r] & E_0 \arrow[r] &  \mathcal{A}_U\otimes N\F  \arrow[r] & 0.
\end{tikzcd}
$$

By definition, the characteristic class $\phi(N\F)$ of $N\F$ on $U$ is the characteristic class $\phi(\xi)$ of the virtual bundle $\xi=\sum_{i=0}^r(-1)^iE_i$.

We set $U_0=U\setminus S$ and $U_1$ a neighbourhood of $S$ in $U$ and consider the covering 
${\mathcal U}=\{U_0,U_1\}$ of $U$. 
On $U_0$, we have the exact sequence of vector bundles
$$
\begin{tikzcd}
0 \arrow[r] & E_r \arrow[r] & E_{r-1} \arrow[r] & \cdots \arrow[r] & E_0 \arrow[r] &  NF  \arrow[r] & 0.
\end{tikzcd}
$$
Let $\nabla$ be a basic connection for $NF$ on $U_0$. We take a connection $\nabla_0^{(i)}$
for each $E_i$ on $U_0$ so that $(\nabla_0^{(r)},\dots,\nabla_0^{(0)},\nabla)$ is compatible with the above sequence (cf.~\cite{BB}). If we denote by $\nabla_0^\bullet$ the family of 
connections $(\nabla_0^{(r)},\dots,\nabla_0^{(0)})$, we have
\begin{equation}\label{eqcomp}
\phi(\nabla_0^\bullet)=\phi(\nabla).
\end{equation}

On $U_1$, we take an arbitrary family $\nabla_1^\bullet=(\nabla_1^{(r)},\dots,\nabla_1^{(0)})$ of connections, each $\nabla_1^{(i)}$ being a connection for $E_i$ on $U_1$. Then the class
$\phi(N\F)=\phi(\xi)$ in $H^{2d}(U;\C)\simeq H^{2d}_D({\mathcal U})$ is represented by the cocycle
\[
(\phi(\nabla_0^\bullet),\phi(\nabla_1^\bullet),\phi(\nabla_0^\bullet,\nabla_1^\bullet)).
\]
If $d\ge n-k+1$, then by \eqref{eqcomp} and the Bott vanishing theorem, we have $\phi(\nabla_0^\bullet)=0$. Thus the above cocycle defines a class in $H^{2d}_D({\mathcal U},U_0)\simeq H^{2d}(U,U\setminus S;\C)$, which is called 
the localization of $\phi(N\F)$ by $\F$ at $S$ and is denoted by $\phi_S(N\F,\F)$. 

Complex analytic subsets are subanalytic. Hence, for a complex analytic subset $S$ of a complex \mfd\ $X$, one may choose a $C^1$-triangulation of $X$ compatible with $S$; see, for example, \cite{Shiota1997}. Thus the dualities and intersection products above may be used combinatorially in the analytic setting.

\begin{defi} The {\it Baum--Bott residue}  of $\F$ with respect to $\phi$ at  $S$, denoted by $\Res_\phi(\F,S)$, is defined to be the image of $\phi_S(N\F,\F)$ by
the Alexander isomorphism
\[
\begin{tikzcd}
H^{2d}(U,U\setminus S;\C) \arrow[r,"\sim"] &H_{2n-2d}^{\rm BM}(S;\C).
\end{tikzcd}
\]
\end{defi}

The Baum--Bott residue is expressed as follows. Let $K_0$ be a triangulation of $U$ compatible with $S$ and let $K$, $K'$ and $K^*$ be as before.
By Theorem~\ref{thAlexC}, we have:
 
\begin{teo}\label{thres}
Let $S\subset \Sing(\F)$ be a connected component and $\phi$  a symmetric homogeneous polynomial of degree $d \geq n-k+1$. Then $\mathrm{Res}_{\phi}(\F, S)$ is represented by the cycle
\[  
\sum_{\mathbb{s} \subset  S} \left(  \int_{\mathbb{s}^*\cap  R_1} \phi(\nabla_1^{\bullet})+\int_{\mathbb{s}^*\cap R_{01}} \phi(\nabla_0^{\bullet},\nabla_1^{\bullet})\right) \mathbb{s}, 
\]    
where  $\mathbb{s}$ runs through all the $(2n -2d)$-simplices of $K$ in $S$.
\end{teo}

In particular, if \( \F \) has dimension one and \( S \) is compact, then
\[
\Res_{c_n}(\F, S) = \text{PH}(\F, S),
\]
where \( \text{PH}(\F, S) \) is the Poincar\'e--Hopf index of \( \F \) at \( S \).

In the global situation, let $(S_\lambda)_\lambda$ be the connected components of $\Sing(\F)$. Assume that each $S_\lambda$ admits a neighbourhood $U_\lambda$ such that the $U_\lambda$'s are mutually disjoint and that $N\F$ admits a locally free resolution by real analytic vector bundles on $U_\lambda$. Then we have the following residue formula:
\[
\sum_\lambda(i_\lambda)_*\Res_\phi(\F,S_\lambda)=[X]\frown\phi(N\F)\qquad\text{in}\ \ H_{2n-2d}^{\rm BM}(X;\C),
\]
where $i_\lambda:S_\lambda\hra X$ denotes the inclusion.

\begin{remark}\label{remBB} {\bf 1.} In the case $S$ is compact, the above residue coincides with the usual Baum--Bott residue. Moreover, if $X$ is compact, the residue formula above reduces to the one in \cite[Theorem 2]{BB}.
\smallskip

\noindent
{\bf 2.} More generally, the same construction defines $\Res_\phi(\F,S)$ for any closed subset $S\subset X$ with the following properties: (1) $S$ is a subcomplex with respect to some triangulation; (2) there is a neighbourhood $U$ of $S$ such that $(U\setminus S)\cap\Sing(\F)=\emptyset$ and $N\F$ admits a locally free resolution by real analytic vector bundles on $U$. Thus the chosen support may contain regular points, provided it contains all singularities in a neighbourhood.
\end{remark}

\section{Proof of Theorem \ref{Main-Theorem}}\label{sec:Main-Theorem}

We first record the two facts needed in the proof.

\begin{lema}\label{lem:normal-sequence}
Let $f:V\to X$ be generically transverse to $\F$, set $\G=f^*\F$, and let $Z\subset V$ be a closed analytic subset of codimension at least two such that $(W\setminus Z)\cap\Sing(\G)=\emptyset$ for some neighbourhood $W$ of $Z$. If $f^*T\F$ and $N\G$ are reflexive on $W$, then, after replacing $W$ by a smaller neighbourhood of $Z$, the canonical sequence on $W\setminus Z$ extends to an exact sequence
\begin{equation}\label{eq:normal-pullback}
0\longrightarrow f^*T\F\longrightarrow f^*TX\overset{\psi}{\longrightarrow}N\G\longrightarrow0.
\end{equation}
In particular, the natural morphism $f^*N\F\to N\G$ is an isomorphism.
\end{lema}

\begin{proof}
We work on $W$ and put $W^0=W\setminus Z$. By \eqref{eq:normal-regular}, there is an isomorphism $\chi:N\G^0\simeq f^*N\F^0$. Since the pull-back of \eqref{eq:foliation} is exact on $W^0$, we have a surjective morphism
\[
\psi^0:=\chi^{-1}\circ f^*\pi:f^*TX^0\longrightarrow N\G^0
\]
with kernel $f^*T\F^0$. 
The sheaf $\mathcal Hom(f^*TX,N\G)$ is reflexive. Since $\cod_W Z\geq2$, the morphism $\psi^0$ extends uniquely to a morphism $\psi:f^*TX\to N\G$. On $W^0$ we have $\psi^0\circ df=\pi_\G$, where $\pi_\G:TW\to N\G$ is the canonical surjection. Hence $\psi\circ df=\pi_\G$ on $W$, and $\psi$ is surjective.

Set $K=\ker\psi$. Since $f^*TX$ is locally free and $N\G$ is reflexive, the depth lemma shows that $K$ is reflexive. Let $a:f^*T\F\to f^*TX$ be the natural morphism. It is injective on $W^0$; hence $\ker a$ is supported on $Z$. Since $f^*T\F$ is reflexive, and therefore torsion-free, $a$ is injective. Moreover, $\psi\circ a=0$, because this equality holds on $W^0$. Thus $f^*T\F$ is a subsheaf of $K$. These two reflexive sheaves agree on $W^0$, and therefore agree on $W$. This proves \eqref{eq:normal-pullback}.

The right exact pull-back of \eqref{eq:foliation} gives
\[
f^*T\F\longrightarrow f^*TX\longrightarrow f^*N\F\longrightarrow0.
\]
Since the first morphism is injective, its cokernel is both $f^*N\F$ and $N\G$. Hence $f^*N\F\simeq N\G$.
\end{proof}

\begin{lema}\label{lem:tor-rigidity}
Let $f:V\to U$ be a real analytic map between real analytic manifolds, and let $\mathcal M$ be a coherent $\mathcal A_U$-module. If
$
\operatorname{Tor}^{f^{-1}\mathcal A_U}_1(\mathcal A_V,f^{-1}\mathcal M)=0,
$
then
\[
\operatorname{Tor}^{f^{-1}\mathcal A_U}_j(\mathcal A_V,f^{-1}\mathcal M)=0
\qquad\text{for all }j\geq1.
\]
\end{lema}

\begin{proof}
The assertion is local. Let $y\in V$, $x=f(y)$, and factor $f$ by its graph,
\[
V\overset{\Gamma_f}{\longrightarrow}U\times V\overset{p}{\longrightarrow}U.
\]
Set $R=\mathcal A_{U,x}$, $C=\mathcal A_{U\times V,(x,y)}$ and $B=\mathcal A_{V,y}$. The homomorphism $R\to C$ is flat, $C$ is a regular local ring, and $B$ is a finite $C$-module. If $M=\mathcal M_x$, flat base change gives
\begin{equation}\label{eq:tor-base-change}
\operatorname{Tor}^R_j(B,M)\simeq\operatorname{Tor}^C_j(B,C\otimes_RM).
\end{equation}
Both $B$ and $C\otimes_RM$ are finite $C$-modules. The rigidity theorem for Tor over regular local rings \cite{Lichtenbaum} says that the vanishing of $\operatorname{Tor}^C_i$ implies the vanishing of $\operatorname{Tor}^C_j$ for every $j\geq i$. Taking $i=1$ and using \eqref{eq:tor-base-change} proves the lemma.
\end{proof}

\begin{proof}[Proof of Theorem~\ref{Main-Theorem}]
Since the assertion is local along $S$ and $Z$, we replace $X$ and $V$ by suitable neighbourhoods and assume that $\Sing(\F)=S$, $\Sing(\G)\subset Z$ and $f(V)\subset U$. By condition \rm{(i)}, choose a real analytic locally free resolution
\[
0\longrightarrow E_r\longrightarrow\cdots\longrightarrow E_1\overset{\rho_1}{\longrightarrow}\mathcal A_U\otimes T\F\longrightarrow0.
\]
Since $\mathcal A_U\otimes_{\mathcal O_U}-$ is exact on coherent holomorphic sheaves \cite[Proposition~6.1]{BB}, this gives a locally free resolution
\begin{equation}\label{ch4exact}
0\longrightarrow E_r\longrightarrow\cdots\longrightarrow E_1\overset{\rho'_1}{\longrightarrow}\mathcal A_U\otimes TX\longrightarrow\mathcal A_U\otimes N\F\longrightarrow0,
\end{equation}
where $\rho'_1=(1\otimes\iota)\circ\rho_1$.

For an $\mathcal A_U$-module $\mathcal T$ we write
$
f^*\mathcal T=\mathcal A_V\otimes_{f^{-1}\mathcal A_U}f^{-1}\mathcal T.
$
For a coherent $\mathcal O_U$-module $\mathcal S$ we have
\begin{equation}\label{eq:analytic-hol-pullback}
f^*(\mathcal A_U\otimes_{\mathcal O_U}\mathcal S)
=\mathcal A_V\otimes_{\mathcal O_V}f^*\mathcal S.
\end{equation}
By condition \rm{(iii)},
\[
0\longrightarrow f^*T\F\longrightarrow f^*TX\longrightarrow N\G\longrightarrow0
\]
is exact. Since holomorphic pull-back is right exact, comparison with the pull-back of \eqref{eq:foliation} also gives a natural isomorphism
\begin{equation}\label{eq:normal-identification}
f^*N\F\simeq N\G.
\end{equation}
Tensoring the exact sequence in condition \rm{(iii)} with $\mathcal A_V$ over $\mathcal O_V$ remains exact by \cite[Proposition~6.1]{BB}. Apply the real analytic pull-back to
\[
0\longrightarrow\mathcal A_U\otimes T\F\longrightarrow\mathcal A_U\otimes TX\longrightarrow\mathcal A_U\otimes N\F\longrightarrow0.
\]
The beginning of the long exact sequence of Tor, together with \eqref{eq:analytic-hol-pullback}, gives
\[
0\longrightarrow
\operatorname{Tor}^{f^{-1}\mathcal A_U}_1
(\mathcal A_V,f^{-1}(\mathcal A_U\otimes N\F))
\longrightarrow
\mathcal A_V\otimes f^*T\F
\longrightarrow
\mathcal A_V\otimes f^*TX.
\]
The last morphism is injective. Hence the first Tor vanishes. Lemma~\ref{lem:tor-rigidity} then gives
\begin{equation}\label{eq:all-tor-zero}
\operatorname{Tor}^{f^{-1}\mathcal A_U}_j
(\mathcal A_V,f^{-1}(\mathcal A_U\otimes N\F))=0
\qquad (j\geq1).
\end{equation}
It follows that the pull-back of \eqref{ch4exact} is exact. Using \eqref{eq:normal-identification}, it takes the form
\begin{equation}\label{ch4exact2}
0\longrightarrow f^*E_r\longrightarrow\cdots\longrightarrow f^*E_1
\longrightarrow\mathcal A_V\otimes f^*TX
\longrightarrow\mathcal A_V\otimes N\G\longrightarrow0.
\end{equation}
Thus \eqref{ch4exact2} is a real analytic locally free resolution of $\mathcal A_V\otimes N\G$.

Put $U_0=U\setminus S$ and choose a neighbourhood $U_1$ of $S$ in $U$. Let $\mathcal U=\{U_0,U_1\}$. As in Subsection~\ref{ssBB}, choose a basic connection $\nabla$ for $N\F$ on $U_0$ and families $\nabla_0^\bullet$ and $\nabla_1^\bullet$ of connections for the vector bundles occurring in \eqref{ch4exact}, compatible with the resolution on $U_0$. The localized class $\phi_S(N\F,\F)$ is represented by
$
(\phi(\nabla_1^\bullet),\phi(\nabla_0^\bullet,\nabla_1^\bullet)).
$

Let $W_i=f^{-1}(U_i)$ and $\mathcal W=\{W_0,W_1\}$. On $W_0$, the isomorphism $N\G^0\simeq f^*N\F^0$ identifies the connection induced by $f^*\nabla$ with a basic connection for $N\G^0$. By \eqref{ch4exact2}, the pulled-back families $f^*\nabla_0^\bullet$ and $f^*\nabla_1^\bullet$ are admissible for the localization of $\phi(N\G)$ by $\G$. Consequently $\phi_Z(N\G,\G)$ is represented by
$
(\phi(f^*\nabla_1^\bullet),\phi(f^*\nabla_0^\bullet,f^*\nabla_1^\bullet)).
$
By naturality of Chern--Weil forms and difference forms,
\begin{equation}\label{eq:localized-pullback}
\phi_Z(N\G,\G)=f^*\phi_S(N\F,\F)
\qquad\text{in }H^{2d}(V,V\setminus Z;\C).
\end{equation}

Finally, the localized intersection product with respect to $f$ is, by definition, the Borel--Moore homological counterpart of relative pull-back under Alexander duality. Hence the diagram
\[
\begin{tikzcd}
H^{2d}(U,U\setminus S;\C)\arrow[r,"A_{U,S}"]\arrow[d,"f^*"'] & H^{\rm BM}_{2n-2d}(S;\C)\arrow[d,"(V\mCdot_f\ )_Z"]\\
H^{2d}(V,V\setminus Z;\C)\arrow[r,"A_{V,Z}"'] & H^{\rm BM}_{2n'-2d}(Z;\C)
\end{tikzcd}
\]
is commutative. Applying it to \eqref{eq:localized-pullback} proves \eqref{formula}.
\end{proof}

\section{Some consequences}\label{sec:consequences}

We record two immediate forms of Theorem~\ref{Main-Theorem}, followed by a push-forward formula for proper generically finite maps.

\begin{cor}\label{cor:inclusion}
Let $i:V\hookrightarrow X$ be the inclusion of a complex submanifold which is generically transverse to $\F$. Let $S$ and $Z=V\cap S$ satisfy the hypotheses of Theorem~\ref{Main-Theorem}. Then
\[
\Res_{\phi}(i^*\F,Z)=\left(V\Cdot\Res_{\phi}(\F,S)\right)_Z.
\]
\end{cor}

\begin{cor}\label{cor:isolated}
In the situation of Corollary~\ref{cor:inclusion}, suppose that $\dim V=n-k+1$ and $Z=\{p\}$. Then $i^*\F$ is one-dimensional near $p$ and
\[
\Res_{\phi}(i^*\F,p)=\left(V\Cdot\Res_{\phi}(\F,S)\right)_p.
\]
For $\phi=c_{n-k+1}$, the residue on the left is the local residue of the vector field defining $i^*\F$; in the isolated case it is represented by the corresponding Grothendieck residue.
\end{cor}

\begin{cor}\label{cor:proper-degree}
Assume the hypotheses of Theorem~\ref{Main-Theorem}, with $\dim V=\dim X=n$, and suppose in addition that $f:V\to X$ is proper, surjective and generically finite of degree $m>0$. Then $f|_Z:Z\to S$ is proper and
\begin{equation}\label{eq:proper-degree}
(f|_Z)_*\Res_{\phi}(f^*\F,Z)=m\,\Res_{\phi}(\F,S)
\in H^{\rm BM}_{2n-2d}(S;\C).
\end{equation}
If $\phi$ has rational coefficients, then
\[
\Res_{\phi}(\F,S)\text{ is rational}
\quad\Longleftrightarrow\quad
\Res_{\phi}(f^*\F,Z)\text{ is rational}.
\]
In particular, rationality is invariant under every proper bimeromorphic map satisfying the hypotheses of Theorem~\ref{Main-Theorem}.
\end{cor}

\begin{proof}
Let $U$ be as in Theorem~\ref{Main-Theorem} and put
\[
\alpha=\phi_S(N\F,\F)\in H^{2d}(U,U\setminus S;\C).
\]
By Theorem~\ref{Main-Theorem}, the localized class of $f^*\F$ on $f^{-1}(U)$ is $f^*\alpha$. Under Alexander duality the corresponding residues are represented by cap products with the Borel--Moore fundamental classes of the ambient manifolds. Since the restriction $f:f^{-1}(U)\to U$ is proper, the projection formula gives
\[
(f|_Z)_*\bigl([f^{-1}(U)]\frown f^*\alpha\bigr)
=f_*[f^{-1}(U)]\frown\alpha.
\]
Since $f$ is generically finite of degree $m$, one has $f_*[f^{-1}(U)]=m[U]$ in Borel--Moore homology. Hence
\[
(f|_Z)_*A_{f^{-1}(U),Z}(f^*\alpha)=mA_{U,S}(\alpha),
\]
which is \eqref{eq:proper-degree}. Proper push-forward is defined with rational coefficients, so rationality of the pulled-back residue implies rationality of $m\Res_{\phi}(\F,S)$, and therefore of $\Res_{\phi}(\F,S)$. The converse is Corollary~\ref{cor:rationality-pullback}.
\end{proof}

\begin{exe}\label{ex:finite-etale-rationality}
Put $q=n-k$ and $r=q+1=n-k+1$. Let $Z$ be a compact connected complex manifold of dimension $k-1=n-r$, let $L$ be a holomorphic line bundle on $Z$, and choose non-zero integers $a_1,\dots,a_r$. For $a<0$ we use the convention $L^{\otimes a}=(L^\vee)^{\otimes(-a)}$. Set
\[
E=L^{\otimes a_1}\oplus\cdots\oplus L^{\otimes a_r},
\quad
X=\operatorname{Tot}(E),
\]
and identify $Z$ with the zero section. On a local trivialization of $E$, with coordinates $(z_1,\dots,z_{k-1},w_1,\dots,w_r)$, consider the subsheaf generated by
\[
\frac{\partial}{\partial z_1},\dots,
\frac{\partial}{\partial z_{k-1}},
\qquad
Y=\sum_{\nu=1}^{r}a_\nu w_\nu\frac{\partial}{\partial w_\nu}.
\]
As in \cite[Example~4]{BB}, these local subsheaves glue to a full integrable subsheaf $T\F\subset TX$, and
$
\operatorname{codim}\F=q,
\ 
\Sing(\F)=Z.
$
Let $\phi$ be symmetric and homogeneous of degree $d$, with
$
q+1=r<d<n,
$
and put $x=c_1(L)$. Baum and Bott compute the residue in this case as the class whose Poincar\'e dual on $Z$ is \cite[eq.~(11.27)]{BB}
\begin{equation}\label{eq:BB-example4-general}
\operatorname{PD}_{Z}\Res_{\phi}(\F,Z)
=
\frac{\phi(a_1,\dots,a_r,0,\dots,0)}{a_1\cdots a_r}\,
x^{d-r}.
\end{equation}
If $\phi$ has rational coefficients, then $\Res_{\phi}(\F,Z)$ is rational. If, in addition, the coefficient in \eqref{eq:BB-example4-general} is non-zero and $x^{d-r}\neq0$, the residue is non-zero.

Now let,
$
g:Z'\longrightarrow Z
$
be a finite \'etale covering of degree $\delta$, set $L'=g^*L$ and
\[
E'=g^*E=(L')^{\otimes a_1}\oplus\cdots\oplus(L')^{\otimes a_r},
\qquad
X'=\operatorname{Tot}(E').
\]
There is a Cartesian diagram
\[
\begin{tikzcd}
X' \arrow[r,"f"] \arrow[d] & X \arrow[d]\\
Z' \arrow[r,"g"'] & Z,
\end{tikzcd}
\]
in which $f$ is finite \'etale of degree $\delta$. Since $f$ is \'etale, the compatibility of the tangent and normal sequences required in Theorem~\ref{Main-Theorem} is automatic, and its remaining local hypotheses are satisfied in this model. The pulled-back foliation $\G=f^*\F$ is the same Baum--Bott model over $Z'$, with singular set the zero section $Z'$. Hence, writing $x'=c_1(L')=g^*x$, formula \eqref{eq:BB-example4-general} gives
\[
\operatorname{PD}_{Z'}\Res_{\phi}(\G,Z')
=
\frac{\phi(a_1,\dots,a_r,0,\dots,0)}{a_1\cdots a_r}\,(x')^{d-r}
=
g^*\!\left(\operatorname{PD}_{Z}\Res_{\phi}(\F,Z)\right).
\]
Since $f|_{Z'}=g$, the projection formula gives
\[
(f|_{Z'})_*\Res_{\phi}(\G,Z')
=\delta\,\Res_{\phi}(\F,Z),
\]
This agrees with Corollary~\ref{cor:proper-degree}. In this family the functoriality formula and the equivalence of rationality before and after pull-back are both read directly from the Baum--Bott formula.

As a concrete specialisation, let $Z=A$ be an abelian variety of dimension $k-1$, let $L$ be symmetric, and take $g=[N]:A\to A$. Then $\deg[N]=N^{2(k-1)}$ and
\[
c_1([N]^*L)^{d-r}=N^{2(d-r)}c_1(L)^{d-r},
\]
and, identifying the pulled-back zero section with $A$, Corollary~\ref{cor:proper-degree} becomes
\[
[N]_*\Res_{\phi}(f^*\F,A)
=N^{2(k-1)}\Res_{\phi}(\F,A).
\]
The scaling on the Poincar\'e-dual side and the degree factor in the push-forward formula are therefore explicit.
\end{exe}

Formula \eqref{BBmax} is the classical transverse case. Theorem~\ref{Main-Theorem} does not require the expected-dimension condition on $S$ or compactness; its essential hypothesis is the exact compatibility of the normal sequence under pull-back. Lemma~\ref{lem:normal-sequence} gives a sufficient reflexivity criterion for this compatibility.

\end{document}